\documentclass{amsart}
\usepackage{amsmath}
\usepackage{amscd}
\usepackage{amssymb}

\newtheorem{theorem}{Theorem}[section]
\newtheorem{theorem/definition}{Theorem/Definition}[section]

\newtheorem{lemma}{Lemma}[section]
\newtheorem{corollary}{Corollary}[section]
\theoremstyle{remark}
\newtheorem{remark}{Remark}[section]

\theoremstyle{definition}

\begin{document}
\title
{A note on compact K\"{a}hler-Ricci flow with positive bisectional
curvature\footnote{}}
\author{Huai-Dong Cao and Meng Zhu}
\address{Department of Mathematics\\ Lehigh University\\
Bethlehem, PA 18015} \email{huc2@lehigh.edu \& mez206@lehigh.edu}

\date{}

\begin{abstract}
 We show that for any solution $g_{i\bar j}(t)$ to the K\"ahler-Ricci
flow with positive bisectional curvature $R_{i\bar i j\bar
j}(t)>0$ on a compact K\"ahler manifold $M^n$, the bisectional
curvature has a uniform positive lower bound $R_{i\bar i j\bar
j}(t)>C>0$. As a consequence, $g_{i\bar j}(t)$ converges
exponentially fast in $C^{\infty}$ to an K\"ahler-Einstein metric
with positive bisectional curvature as $t\to \infty$, provided we
assume the Futaki-invariant of $M^n$ is zero. This improves a
result of D. Phong, J. Song, J. Sturm and B. Weinkove \cite{PSSW}
in which they assumed the stronger condition that Mabuchi K-energy
is bounded from below.
\end{abstract}
\maketitle

\footnotetext[1]{Research partially supported by NSF grants
DMS-0354621 and DMS-0506084.}

\section{The Results}

\hspace{.75cm}

We consider solutions $g_{i\bar j} (t)$ to the normalized
K\"ahler-Ricci flow
$$
\frac{\partial g_{i\bar j}}{\partial t}=-R_{i\bar j}+ g_{i\bar j}
\eqno (1.1)
$$
with positive holomorphic bisectional curvature $R_{i\bar i j\bar
j}(t)>0$ on a compact K\"ahler manifold $M^n$. Our main result is
a uniform lower bound on the bisectional curvature of $g_{i\bar j}
(t)$.

\begin{theorem}

Let $g_{i\bar j}(t), 0\le t<\infty,$ be a solution to the
normalized K\"ahler-Ricci flow (1.1) with positive holomorphic
bisectional curvature $R_{i\bar i j\bar j}(t)>0$ on a compact
K\"ahler manifold $M^n$. Then, there exists a positive constant
$C_1>0$ such that
$$R_{i\bar i j\bar j}(t)>C_1$$
for all $0\le t<\infty$.
\end{theorem}

In the Riemannian case, powerful curvature pinching estimates were
first proved by Hamilton for Ricci flow on compact $3$-manifolds
with positive Ricci curvature \cite{Ha82} and on $4$-manifolds
with positive curvature operator \cite{Ha86}. More recently,
curvature pinching estimates were obtained by B\"ohm-Wilking
\cite{BW} for Ricci flow on compact $n$-manifolds ($n\ge 5$) with
positive curvature operator (see also Brendle-Schoen \cite{BS}
with a certain positive curvature condition that is weaker than
$1/4$-pinch). These pinching results are proved by using
Hamilton's advanced maximum principle, and a key step is to show
that the pinching holds under the reaction ODE system of the
evolution equation of the curvature tensor. However, the K\"ahler
case appears to be more subtle since the reaction ODE system of
the curvature tensor doesn't even preserve the lower bound of the
bisectional curvature, even in complex dimension $2$. In fact,
Hamilton \cite{Ha95F} showed that under the reaction ODE system,
there could be a solution that emerges from the curvature operator
matrix of $\Bbb{CP}^2$ and approaches the curvature matrix of
$\Bbb{C}P^{1}\times \Bbb{C}$. Of course this doesn't mean it
actually happens to the curvature evolution PDE system, but a
different argument from the Riemannian case is needed. Our proof
of the uniform lower bound for the bisectional curvature relies on
a result of H.-L. Gu \cite{Gu} (see Lemma 2.1 below) which in turn
depends on the strong maximum principle developed by
Brendle-Schoen \cite{BS2} for certain degenerate parabolic
operators. We also use a lemma from Chen-Sun-Tian \cite{CST}.
While we are able to prove a uniform positive lower bound on the
bisectional curvature, it remains very interesting to see if one
can show the pinching of holomorphic
sectional or bisectional curvature.\\

Based on Theorem 1.1 and the work of \cite{PSSW}, we obtain the
following

\begin{theorem}
Let $(M^{n}, g_{0})$ be a compact K\"{a}hler manifold with
positive bisectional curvature. Assume that the Futaki invariant
of $M^n$ is zero. Then the solution $g_{i\bar{j}}(t)$ to (1.1)
with initial metric $g_{0}$ converges in $C^{\infty}$-norm
exponentially fast to a K\"{a}hler-Einstein metric $\tilde{g}$. In
particular, $M^{n}$ is biholomorphic to the complex projective
space $\Bbb{C}P^n$ and $\tilde{g}$ is the Fubini-Study metric.
\end{theorem}

Theorem 1.2 improves one of the main results of \cite{PSSW}, in
which it is assumed that the Mabuchi $K$-energy is bounded below.
As pointed out in \cite{PSSW}, the vanishing of the Futaki
invariant is, at least a priori, much weaker than the K-energy
being bounded below. Of course, one would like to drop the
assumption on Futaki invariant as well.

\begin{remark}
For $n=1$, Hamilton \cite{Ha88} showed that under the Ricci flow,
one can deform any metric of positive Gauss curvature on a compact
surface $\Sigma$ with positive Euler number to a metric of
constant positive Gauss curvature. Later in \cite{Chow}, Ben Chow
showed that any initial metric on $\Sigma$ becomes positive curved
after some finite time, hence converges to a metric of constant
positive Gaussian curvature. For $n=2$, Theorem 1.2 was derived in
\cite{PSSW} without using Theorem 1.1 due to the basic fact in
K\"ahler geometry that in complex dimension 2 nonnegative
bisectional curvature implies nonnegative curvature operator on
(2,0)-tensors (also see, e.g., \cite{Z})
\end{remark}

\medskip
{\bf Acknowledgement.} We would like to thank Richard Hamilton,
D.H. Phong, Jacob Sturm, Xiaofeng Sun and Valentino Tosatti for
their interest in our work.

\section{The Proofs}

\noindent {\bf Proof of Theorem 1.1}. Let $g_{i \bar j}(t)$, $0\le
t <\infty$, be a solution to the normalized K\"ahler-Ricci flow
(1.1) with positive bisectional curvature on a compact $M^n$.

\smallskip
{\bf Claim} {\sl There exists some positive constant $C>0$, such
that
$$R_{i \bar j}(t)>Cg_{i \bar j}(t) \eqno(2.1)$$
for all $t\ge 0$.}
\smallskip

We argue by contradiction. Suppose the claim is not true. Then we
can find a sequence of positive number $\epsilon_k \to 0$, and a
sequence of points $\{(x_k,t_k)\}_{k=1}^{\infty}$ in space-time
with $x_k\in M$ and $t_k \to \infty$ as $k \to \infty$ such that
$$\min R_{i\bar{j}}(x_k,t_k)\le \epsilon_k. \eqno(2.2)$$
Now we can choose a unitary frame $E^{k}=\{e^k_1,\cdots, e^k_n \}$
at the point $x_k$ and the time $t_k$ so that
$$R_{1\bar{1}}(x_k,t_k)=\min_{1\le i\le n}
R_{i\bar{i}}(x_k,t_k).$$ By the work of B.-L. Chen, X.-P. Zhu and
the first author \cite{CCZ} (or by Perelman, see \cite{ST}), we
know that the diameter and the curvature tensor of $g_{i\bar
j}(t)$ are uniformly bounded. Also injectivity radius is uniformly
bounded by Perelman's noncollapsing theorem \cite{P}. Thus we can
apply Hamilton's compactness theorem in \cite{Ha95} so that a
subsequence of the compact marked solutions $\{(M^n, J,
g_{i\bar{j}}(t_k+t),x_k, E^k)\}$ converges in $C^{\infty}$ norm,
in the Cheeger-Gromov sense, to a compact marked solution $(M,
\tilde{J}, \tilde{g}_{i\bar{j}}(t), \tilde{x},\tilde{E})$ to (1.1)
with nonnegative bisectional curvature
$\tilde{R}_{i\bar{i}j\bar{j}}\ge 0$ and
$$\tilde{R}_{1\bar{1}}(\tilde{x}, 0)=\lim_{k\to \infty}R_{1\bar{1}}(x_k, t_k)=0.\eqno(2.3)$$
Here $\tilde{E}$ is a unitary frame at the marked point
$\tilde{x}$ at $t=0$, and $\tilde{J}$ is a complex structure on
$M^n$, possibly different from $J$. Moreover, we know that $(M,
\tilde{J}, \tilde{g}_{i\bar{j}}(t))$ is a gradient shrinking
K\"{a}hler-Ricci soliton.

Now we use the following result of H.-L. Gu \cite{Gu}.

\begin{lemma} Given any K\"ahler metric $h_{i\bar j}$ with nonnegative bisectional
curvature on a compact, irreducible, simply connected K\"ahler
manifold $N^n$. Then, under the normalized K\"ahler-Ricci flow
(1.1), either the bisectional curvature becomes positive
everywhere after a short time, or $(N^n, h_{i\bar j})$ is
isometrically biholomorphic to a Hermitian symmetric space of rank
$\ge 2$.
\end{lemma}

We know that our manifold $M^n$ is simply connected. Also, by a
theorem of Bishop and Goldberg \cite{BG} (see also Theorem 4,
\cite{BG2}), $b_2(M^n)=1$. Hence, $M^n$ is irreducible. On the
other hand, since $(M, \tilde{g})$ is a shrinking soliton with
nonnegative bisectional curvature, (2.3) tells us that for each
$t\in [0, \infty)$, the bisectional curvature $\tilde{R}_{i\bar i
j\bar j}(t)$ vanishes somewhere at some point on $M^n$ at $t$.
Therefore, Lemma 2.1 implies that $(M^n, \tilde{J}, \tilde{g})$ is
isometrically biholomorphic to a Hermitian symmetric space of rank
$\geq 2$. This would lead to a contradiction, either by Mok's
solution to the Generalized Frenkel conjecture \cite{Mok}, or
directly because $\tilde{g}$ being K\"ahler-Einstein, i.e.,
$\tilde{R}_{i\bar j}=\tilde{g}_{i\bar j},$ in turn implies that
$$||R_{i\bar{j}}-g_{i\bar{j}}||_{C^0}(t_k)\to 0,$$ contradicting
(2.2). This finishes the proof of the Claim.

On the other hand, by a direct computation of the evolution of
$$S_{i\bar jk\bar l}(t)=R_{i\bar{j}k\bar{l}}(t)-c\big[g_{i\bar{j}}R_{k\bar{l}}
+R_{i\bar{j}}g_{k\bar{l}}+g_{i\bar{l}}R_{k\bar{j}}+R_{i\bar{l}}g_{k\bar{j}}\big](t)$$
for some sufficiently small positive constant $c>0$ and applying
the maximum principle, Chen-Sun-Tian \cite{CST} (see Lemma 6 in
\cite{CST}) showed that if the Ricci curvature of
$g_{i\bar{j}}(t)$ has a uniform positive lower bound then the
bisectional curvature of
$g_{i\bar{j}}(t)$ has a uniform positive lower bound. Thus Theorem 1.1 follows. \\

In \cite{PS}, Phong and Sturm first considered the smallest
positive eigenvalue $\lambda$ of the operator $L=-g^{i\bar
j}\nabla_{i}\nabla_{\bar j}$ acting on smooth (1,0) vector fields.
It has been demonstrated in \cite{PS} and \cite{PSSW, PSSW2} that
a uniform positive lower bound on $\lambda$ plays an important
role in showing the convergence of the K\"ahler-Ricci flow.

Combining Theorem 1.1 above with Lemma 6 and Lemma 2 in
\cite{PSSW}, we derive a uniform positive lower bound on $\lambda$
for solutions to the K\"ahler-Ricci flow (1.1) with positive
bisectional curvature.

\begin{corollary} Under the same assumption as in Theorem 1.1, there
exists a positive constant $C_2>0$ such that the smallest positive
eigenvalue $\lambda$ satisfies
$$\lambda\ge C_2$$ for all $t\ge 0$.
\end{corollary}

\begin{remark} Under the additional assumption that the Mabuchi
K-energy is bounded from below, both Theorem 1.1 and Corollary 2.1
were proved in \cite{PSSW}.
\end{remark}

Now we are ready to prove Theorem 1.2. \\

\noindent \textbf{Proof of Theorem 1.2}. Let $(M^{n}, g_{0})$ be a
compact K\"{a}hler manifold with positive bisectional curvature.
Then, by the work of the first author \cite{Cao}, we know that the
solution $g_{i\bar j}(t)$ to the normalized K\"ahler-Ricci flow
(1.1) with initial metric $g_{0}$ exists for all time $0\le
t<\infty$. Moreover, it follows from Bando \cite {Ba} (for $n=3$)
and Mok \cite{Mok} (for $n\ge 3$) that $g_{i\bar j}(t)$ has
positive bisectional curvature for all $t\ge 0$.

We will need the following result from \cite{PSSW}

\begin{lemma}{\bf (Phong-Song-Sturm-Weinkove \cite{PSSW})}
Suppose we have a nonsingular solution to the
normalized K\"ahler-Ricci flow (1.1) on a compact K\"ahler
manifold with positive first Chern class. Assume that the Futaki
invariant is zero and the smallest positive eigenvalue $\lambda$
of the operator $L$ has a uniform positive lower bound. Then the
solution converges to a K\"ahler-Einstein metric exponentially
fast in $C^{\infty}$-norm.
\end{lemma}

Now our solution $g_{i \bar j}(t)$ is nonsingular by Cao-Chen-Zhu
\cite{CCZ} (or Perelman, see \cite{ST}), and Futaki invariant is
zero by assumption. Moreover, by Corollary 2.1, the smallest
positive eigenvalue $\lambda$ of the operator $L=-g^{i\bar
j}\nabla_{i}\nabla_{\bar j}$ has uniform positive lower bound.
Thus it follows from Lemma 2.2 that $g_{i\bar{j}}(t)$ converges
exponentially fast in $C^{\infty}$ to a K\"{a}hler-Einstein metric
$\tilde{g}$. Also, by Theorem 1.1, $\tilde{g}$ has positive
bisectional curvature. Therefore, by a theorem of M. Berger
\cite{Be} (see also Theorem 5, \cite{BG2}), $\tilde{g}$ is of
constant holomorphice sectional curvature. The proof of Theorem 2
is completed.

\end{document}